\documentclass[12pt]{article}
\usepackage{amsmath,amsthm,amssymb, amstext, amsopn, amsxtra}
\usepackage[all]{xy}
\usepackage{latexsym}
\usepackage{amsfonts}
\begin{document}
\newcommand{\nt}{\noindent}
\newcommand{\bs}{\bigskip}
\newcommand{\ms}{\medskip}
\newcommand{\mk}{\medskip}
\newcommand{\sk}{\smallskip}
\newcommand{\ep}{\varepsilon}
\newcommand{\m}{{\mathfrak m}}
\newcommand{\p}{{\mathfrak p}}
\newcommand{\dd}{\delta}
\newcommand{\A}{{\mathbb A}}
\newcommand{\R}{{\mathcal R}}
\newcommand{\M}{{\mathbb M}}
\newcommand{\Rstar}{$\,(R,\,{}^*\,)\,$}
\newcommand{\starring}{$\,{\ast}\,$-ring}
\newcommand{\starrings}{$\,{\ast}\,$-rings}
\newcommand{\ap}{\alpha^{\,\prime}}
\newcommand{\bp}{\beta^{\,\prime}}
\newcommand{\isom}{\cong}

\bs

\sk
\begin{center}
\Large{\bf A New Determinantal Formula for Three Matrices}

\bigskip
\large{Dinesh Khurana and T.\,Y.~Lam} 
\end{center}

\mk
\begin{abstract}
\begin{small}
\nt
For any three $\,n\times n\,$ matrices $\,A,B,X\,$ over a commutative
ring $\,S$, we prove that $\,{\rm det}\,(A+B-AXB)={\rm det}\,(A+B-BXA)
\in S$. This apparently new formula may be regarded as a ``ternary
generalization'' of Sylvester's classical determinantal formula
$\,{\rm det}\,(I_n-AB)={\rm det}\,(I_n-BA)\,$ for any pair of
$\,n\times n\,$ matrices $\,A,B\,$ over $\,S$. The discovery and proof
of this new generalization of Sylvester's formula were prompted by the
authors' past and present work on the element-wise theory of ring
elements of stable range one, in [KL$_1$] and [KL$_2$].
\end{small}
\end{abstract}

\bs\nt 
{\bf \S1. \ Introduction}

\bs\indent
For any two $\,n\times n\,$ matrices $\,A,\,B\,$ over a commutative ring
$\,S$, it is well known that $\,{\rm det}\,(I-AB)={\rm det}\,(I-BA)$.
In the standard linear algebra literature, this nice formula is often
cited as ``Sylvester's Determinant Identity'', since it was first stated
by J.\,J.~Sylvester in 1857 (for square matrices over the complex numbers).
In a number of other references, the same formula was sometimes also
called the ``Weinstein-Aronszajn Determinant Identity''; see, for instance,
Tao's expository survey [Ta], or Penn's Youtube podcast [Pn].  In our
recent work on the relatively new theory of ring elements of stable range
one, we discovered that it is actually possible to prove a generalization
of Sylvester's classical formula; namely, if $\,A,B,X\,$ are three
arbitrary $\,n\times n\,$ matrices over a commutative ring $\,S$, then
$$
{\rm det}\,(A+B-AXB)={\rm det}\,(A+B-BXA)\in S.  \leqno (1.1)
$$
The detailed proof of this new determinantal identity (with accompanying
motivational remarks) will be given in \S2. As is to be totally expected,
(1.1) is a ``ternary generalization'' of Sylvester's classical determinantal
identity, since the latter can be easily retrieved from (1.1) by setting
$\,X=I\,$ (the $\,n\times n\,$ identity matrix) and replacing $\,A,\,B\,$
respectively by $\,I-A\,$ and $\,I-B$.  (Curiously enough, setting
$\,B=I\,$ and $\,X=I+B\,$ would have also worked.) Some easy variations
and applications of the new identity (1.1) will be presented subsequently
in \S3.  While the larger part of this paper was motivated by the authors'
recent work [KL$_2$] on the theory of ring elements of stable range one,
the writing of the present paper is designed to be almost completely
independent of that of [KL$_2$].

\bs
Throughout the paper, the notation $\,{\mathbb M}_n(S)\,$ will be used
to denote the ring of $\,n\times n\,$ matrices over a ring $\,S$, and the
group of invertible $\,n\times n\,$ matrices over $\,S\,$ will be denoted
by $\,{\rm GL}_n(S)$.  As is well known, in the case where $\,S\,$ is a
{\it commutative\/} ring, a matrix $\,A\in {\mathbb M}_n(S)\,$ belongs
to $\,{\rm GL}_n(S)\,$ if and only if $\,{\rm det}\,(A)\,$ is a unit
of $\,S$.  In \S2, the notation $\,E_{ij}\,$ will be used to denote the
(so-called) {\it matrix units\/} in the matrix ring $\,{\mathbb M}_n(S)$.
Other standard terminology and conventions in matrix theory and ring
theory follow largely those in [Bh], [HJ$_1$, HJ$_2$] and [La]. For a
thorough survey on the results and historical developments in the area
of determinantal identities for square matrices, we recommend the
excellent article [BrS] of Brualdi and Schneider.

\mk

\bs\nt
{\bf \S2. \ A Ternary Determinantal Formula}

\bs
Before we proceed to the main body of this paper, some words of
explanation and motivation should be helpful.  In noncommutative
ring theory, there is a delightful result, popularly known as
{\bf Jacobson's Lemma\/}, which states that, for any two elements
$\,a,\,b\,$ in a ring $\,R$, one has $\,1-ab\in {\rm U}(R)\,$
if and only if $\,1-ba \in {\rm U}(R)$, where $\,{\rm U}(R)\,$
denotes the group of units in $\,R$. (This classical lemma should
be best thought of as a purely ring-theoretic result that is
``inspired by'' Sylvester's Determinant Identity.) In the authors'
recent work [KL$_2$] on the theory of ring elements of stable range
one, it was discovered that there is a {\it ternary generalization\/}
of Jacobson's Lemma (fondly called {\bf Super Jacobson's Lemma\/}
in [KL$_2$]) to the effect that, for any three ring elements $\,a,b,x
\in R$, $\,a+b-axb\in{\rm U}(R)\,$ if and only if $\,a+b-bxa\in
{\rm U}(R)$.  Using this new ternary lemma, the authors were able
to prove in [KL$_2$: Theorem 3.1] that an arbitrary ring element
has left stable range one if and only if it has right stable
range one.  While this new left-right symmetry result for ring
elements of stable range one does not directly say anything offhand
about the determinants of matrices, the authors realized for the
first time from Super Jacobson's Lemma that, for three $\,n\times n\,$
matrices $\,A,B,X\,$ over a {\it commutative\/} ring,
$\,{\rm det}\,(A+B-AXB)\,$ must be somehow ``intimately related''
to $\,{\rm det}\,(A+B-BXA)$. After some nontrivial work with block
elementary transformation of matrices, we managed to prove the
remarkable determinantal identity (2.2) in the theorem below, as
a ``ternary generalization'' of Sylvester's classical determinantal
identity for two matrices. Surprisingly to us, a quick search through
the standard literature in matrix theory and determinant theory
(e.g.~[Mu], [Bh], [BrS], and [HJ$_1$, HJ$_2$]) did not turn up
the following result.

\bs\nt
{\bf Theorem 2.1.} {\it For any $\,A,B,X\in R={\mathbb M}_n(S)\,$
over a commutative ring $\,S$, we have}
$$
{\rm det}\,(A+B-AXB)={\rm det}\,(A+B-BXA). \leqno (2.2)
$$
{\it Consequently, $\,A+B-AXB \in {\rm GL}_n(S)\,$ if and only
if $\,A+B-BXA \in {\rm GL}_n(S)$.  However, we may not have
$\,{\rm tr}\,(A+B-AXB)={\rm tr}\,(A+B-BXA)$. Instead, we have
always}
$$
{\rm tr}\,(A+B-AXB)={\rm tr}\,(A+B-XBA)={\rm tr}\,(A+B-BAX),
\leqno (2.3)
$$
{\it while the determinants of the three matrices in $(2.3)$
may be all different.}

\bs\nt
{\bf Proof.} The statement about membership in $\,{\rm GL}_n(S)\,$
follows quickly from (2.2) since a matrix $\,Y\in R\,$ is invertible
in $\,R\,$ iff $\;{\rm det}\,(Y)\in {\rm U}(S)\,$ (the unit group
of the ring $\,S$).  To prove (2.2), let $\,I=I_n\in R\,$ and write
down the following easy matrix identity (which has appeared before,
for instance, in [LN: \S4]):
$$
\begin{pmatrix}I&-B\\I-AX&A\end{pmatrix}
\begin{pmatrix}I&0\\X&I\end{pmatrix}
= \begin{pmatrix}I-BX&-B\\I&A\end{pmatrix}. \leqno (2.4)
$$
Next, we notice that the first and the third block matrices above can
be changed by elementary ``block column transformations'' as follows\,:  
$$
\begin{pmatrix}I&-B\\I-AX&A\end{pmatrix}
\begin{pmatrix}I&B\\0&I\end{pmatrix}
= \begin{pmatrix}I&0\\I-AX&A+B-AXB\end{pmatrix}. \leqno (2.5)
$$
$$
\begin{pmatrix}I-BX & -B\\I & A\end{pmatrix}
\begin{pmatrix}I & -A\\0 & I\end{pmatrix}
= \begin{pmatrix}I-BX & -A-B+BXA\\I&0\end{pmatrix}. \leqno (2.6)
$$
From (2.4), we see that the two left factors on the LHS of (2.5) and
(2.6) have the same determinants.  Thus, the two matrices on the RHS
of (2.5) and (2.6) also have the same determinants.  We are done by
finally noting that the last two determinants are exactly the LHS
and the RHS of (2.2).

\mk
To see that (2.2) may not hold (for $\,n\geq 2\,$ and $\,S\neq (0)$)
if ``determinant'' is replaced by ``trace'', we take $\,A\,$ and
$\,X\,$ to be respectively the matrix units $\,E_{11}\,$ and
$\,E_{12}$, so that $\,AX=E_{12}\,$ and $\,XA=0$.  Choosing
$\,B=s\,E_{21}\,$ with $\,s\in S\setminus \{0\}$, we have
$\,{\rm tr}\,(BXA)=0$, while $\,{\rm tr}\,(AXB)={\rm tr}\,(s\,E_{11})
=s\neq 0$.  From this, it follows that
$$
{\rm tr}\,(A+B-AXB)\neq {\rm tr}\,(A+B-BXA) \leqno (2.7)
$$
in this case.  Nevertheless, the trace equations in (2.3)
do always hold, in view of the well known fact that
$\,{\rm tr}\,(YZ)={\rm tr}\,(ZY)\,$ for any $\,Y,Z
\in {\mathbb M}_n(S)\,$ over the {\it commutative\/} ring $\,S$.

\mk
Our final job is to show that the {\it determinants\/}
of the three matrices $\,P=A+B-AXB$, $\,H=A+B-XBA\,$ and
$\,K=A+B-BAX\,$ in (2.3) may in fact be all different.
This is easily checked by taking, for instance,
$\,A=\begin{pmatrix}0&1\\0&0\end{pmatrix}$,
$\,B=\begin{pmatrix}0&1\\1&0\end{pmatrix}$, and
$\,X=\begin{pmatrix}0&x\\y&1\end{pmatrix}$. 
For these choices, a simple calculation shows that
$\,{\rm det}\,(P)=y-2$,  $\,{\rm det}\,(H)=x-2$,
and $\,{\rm det}\,(K)=2y-2$, which are, in general,
three different elements in the ring $\,S$. \qed

\bs
After giving the above proof for Theorem 2.1, we should point
out, with a good dose of hindsight, that a quicker proof for the
determinantal identity (2.2) is actually possible. Indeed, as far
as proving (2.2) is concerned, we may ``assume'' that the entries
of the matrices $\,A,\,B\,$ and $\,X\,$ are independent commuting
variables, say over the ring $\,{\mathbb Z}$.  By doing so,
we may replace $\,S\,$ by a polynomial ring $\,\hat{S}\,$ over
$\,{\mathbb Z}\,$ in those commuting variables. Over the quotient
field of $\,\hat{S}$, the ``generic matrices'' $\,A\,$ and $\,B\,$
become both {\it invertible.}  Letting $\,P\!:=A+B-AXB\,$ and
$\,Q\!:=A+B-BXA$, we see that
$$
A^{-1}PB^{-1}=A^{-1}(A+B-AXB)\,B^{-1}=B^{-1}+A^{-1}-X, \leqno (2.8)
$$
and similarly, $\,B^{-1}QA^{-1}=B^{-1}(A+B-BXA)\,A^{-1}=B^{-1}+A^{-1}-X$.
Thus, $\,A^{-1}PB^{-1}=B^{-1}QA^{-1}$.  Taking the determinants on
both sides, we may cancel the factors $\,{\rm det}\,(A^{-1})\,$ and
$\,{\rm det}\,(B^{-1})\,$ to conclude that $\,{\rm det}\,(P)
={\rm det}\,(Q)$, thus proving (2.2).

\bs
Recall that, in classical linear algebra, two matrices $\,M,N\in
{\mathbb M}_m(S)\,$ are said to be {\it equivalent\/} if there exist
$\,U,V\in {\rm GL}_m(S)\,$ such that $\,N=UM\,V$.  If $\,U,\,V\,$
can both be chosen in $\,{\rm SL}_m(S)$, we will say more precisely
that $\,M,\,N\,$ are SL-equivalent. By checking through the proof
steps of Theorem 2.1 carefully (reviewing, especially, the three matrix
equations (2.4), (2.5) and (2.6)), we can show that, for $\,P=A+B-AXB\,$
and $\,Q=A+B-BXA\,$ above, the two ``suspended matrices''
$\,{\rm diag}\,(P,I_n)\,$ and $\,{\rm diag}\,(Q,I_n)\,$ are
$\,{\rm SL}$-equivalent in $\,{\mathbb M}_{2n}(S)$.  In the special
case where $\,A,\,B\in {\rm GL}_n(S)$, the equation $\,A^{-1}PB^{-1}
=B^{-1}QA^{-1}\,$ in the last paragraph would even show that $\,P\,$
and $\,Q\,$ themselves are equivalent.  However, by working with a
specific example, say over the polynomial ring $\,S={\mathbb Z}\,[x]$,
with $\,X=I_2\,$ and $\,A=\footnotesize{\begin{pmatrix}0&0\\0&x
\end{pmatrix}}$, $B=\footnotesize{\begin{pmatrix}0&2\\0&x\end{pmatrix}}$,
one can show that $\,P=\footnotesize{\begin{pmatrix}0&2\\0 & 2x-x^2
\end{pmatrix}}\,$ and $\,Q=\footnotesize{\begin{pmatrix}0&2-2x\\0&2x-x^2
\end{pmatrix}}$ are {\it not\/} equivalent in $\,{\mathbb M}_2(S)$.
Verification of this non-equivalence will be left to the interested
reader.

\mk

\bs\nt
{\bf \S3. \ Some More Determinantal Identities, and Examples}

\bs
Concerning ``ternary generalizations'' of Sylvester's Identity, we should
point out that other ``more obvious'' (or perhaps ``better looking'')
forms of such generalizations may not hold true at all. For instance,
given three arbitrary matrices $\,A,X,B\in R={\mathbb M}_n(S)\,$ over
a ring $\,S$, one may ask if the condition $\,I-AXB\in{\rm GL}_n(S)\,$
might be equivalent to the condition $\,I-BXA\in {\rm GL}_n(S)$.
In the special case where $\,X=I$, the answer to this question
is ``yes'' according to Jacobson's Lemma.  On the other hand,
if $\,X\neq I$, the following easy example shows that the answer
to the above question is, in general, ``no''.

\mk\nt
{\bf Example 3.1.} Let $\,A=\footnotesize{\begin{pmatrix}
1&0\\0&0\end{pmatrix}}$, $X=\footnotesize{\begin{pmatrix}0&1\\1&0
\end{pmatrix}}$, and $\,B=\footnotesize{\begin{pmatrix}1&1\\0&0
\end{pmatrix}}$ over any ring $\,S\neq \{0\}$.  Here, $\,I-AXB=I
\in {\rm GL}_2(S)\,$ has determinant $\,1$, but
$\,I-BXA=\footnotesize{\begin{pmatrix}0&0\\0&1\end{pmatrix}}
\notin {\rm GL}_2(S)\,$ has determinant $\,0$. On the other hand,
for the same choices of $\,A,B,X\in {\mathbb M}_2(S)$, both
$\,A+B-AXB=\footnotesize{\begin{pmatrix}2&1\\0&0\end{pmatrix}}$
and $\,A+B-BXA=\footnotesize{\begin{pmatrix}1&1\\0&0\end{pmatrix}}$
have determinant zero, thus reaffirming the truth of Theorem 2.1
in this particular case.

\bs
To further illustrate the meaning and significance of the
``Super Jacobson's Lemma'' mentioned in \S2, we shall state and
prove {\it four\/} of its interesting ``binary specializations''
below, in terms of the determinants and traces of square matrices.

\bs\nt
{\bf Theorem 3.2.} {\it For any $\,A,\,X \in R={\mathbb M}_n(S)\,$
over a commutative ring $\,S$, the following four matrices
in $\,R\,$ have the same determinants as well as the same traces\,}:

\mk\nt
(1) $\,M_1=I-AX+AXA$.\;\;\;        (2) $\,M_2=I-XA+AXA$.\\
(3) $\,M_3=I-AX+A^2 X$.\;\;\;\;\,  (4) $\,M_4=I-XA+X A^2$.

\mk\nt
{\it In particular, in the special case where $\,n\leq 2$,
all four matrices above have the same characteristic polynomials.}

\bs\nt
{\bf Proof.}  First, $\,{\rm det}\,(M_1)={\rm det}\,(M_2)\,$
follows from the equation (2.2) by taking $\,B=I-A$.  Next,
$\,{\rm det}\,(M_1)={\rm det}\,(M_4)\,$ follows from Sylvester's
Identity by writing
$$
\mbox{$M_1=I-A\,(X-XA)$,\,\;\;and\,\;\;$\,M_4=I-(X-XA)\,A$.}
$$
Finally, $\,{\rm det}\,(M_2)={\rm det}\,(M_3)\,$ follows from
Sylvester's Identity again by writing
$$
\mbox{$\,M_2=I-(X-AX)\,A$,\,\;\;and\,\;\;$\,M_3=I-A\,(X-AX)$.}
$$
As for traces, the trace of $\,M_1=I-AX(I-A)\,$ is equal to
the trace of $\,I-(I-A)AX=M_3$, as well as to the trace of
$\,I-X(I-A)A=M_4$. Finally, the trace of $\,M_4=I-XA\,(I-A)\,$
is equal to the trace of $\,I-(I-A)XA=M_2$, so all four matrices
$\,M_i\,$ have the same traces too. The last conclusion of the
theorem now follows trivially from the fact that, in the special
cases $\,n=1\,$ and $\,n=2$, the characteristic polynomial of
an $\,n\times n\,$ matrix $\,M\,$ is completely determined
by $\,{\rm tr}\,(M)\,$ and $\,{\rm det}\,(M)$. \qed

\bs\nt
{\bf Example 3.3.}  To run a random check on the
conclusions of Theorem 3.2, take the two matrices
$\,A=\footnotesize{\begin{pmatrix}1&r\\1&0\end{pmatrix}}$ and
$\,X=\footnotesize{\begin{pmatrix}s&t\\0&0\end{pmatrix}}$
in $\,{\mathbb M}_2(S)\,$ for any $\,r,s,t\,$ in a commutative
ring $\,S$. A quick calculation shows that the four matrices in
Theorem 3.2 are, respectively:
$$
M_1=\footnotesize{\begin{pmatrix}1+t&sr-t\\t&1+sr-t\end{pmatrix}},\;\;
M_2=\footnotesize{\begin{pmatrix}1&0\\s+t&1+sr\end{pmatrix}},\;\;
{\rm and} \;\;
M_3=M_4=\footnotesize{\begin{pmatrix}1+sr&rt\\0&1\end{pmatrix}}.
$$
All four matrices have determinant $\,1+sr\,$ and trace $\,2+sr$,
both of which turned out to be (somehow) independent of the
element $\,t\in S$.

\sk\nt

\bs
\nt Department of Mathematics \\
\nt Panjab University \\
\nt Chandigarh-160014, India

\sk
\nt {\tt dkhurana@pu.ac.in}  

\bs
\nt Department of Mathematics \\
\nt University of California \\
\nt Berkeley, CA 94720, USA   

\sk
\nt {\tt lam@math.berkeley.edu}

\end{document}